\documentclass{ws-ijbc}
\usepackage{ws-rotating}     
\usepackage{graphicx}
\usepackage{graphics}
\usepackage{epstopdf}
\usepackage{amssymb}
\usepackage{amsmath}
\usepackage{subfigure}
\begin{document}

\catchline{}{}{}{}{} 

\markboth{Hui Wang, Xiaoli Chen, Jinqiao Duan}{A Stochastic Pitchfork Bifurcation in Most Probable Phase Portraits}

\title{A Stochastic Pitchfork Bifurcation\\ in Most Probable Phase Portraits }

\author{Hui Wang, Xiaoli Chen}
\address{Center for  Mathematical Sciences  $\&$  School of Mathematics and Statistics\\  $\&$ Hubei Key Laboratory for Engineering Modeling and Scientific Computing,\\ Huazhong University of Science and Technology, Wuhan 430074, China \\
huiwheda@hust.edu.cn,\hspace{1cm} xlchen@hust.edu.cn}

\author{Jinqiao Duan}
\address{Department of Applied Mathematics, Illinois Institute of Technology, Chicago, IL 60616, USA \\$\&$ Center for  Mathematical Sciences, \\
 Huazhong University of Science and Technology, Wuhan 430074, China \\    duan@iit.edu}

\maketitle

\begin{history}
\received{(to be inserted by publisher)}
\end{history}

\begin{abstract}
We study   stochastic   bifurcation for a   system under multiplicative stable L\'evy noise (an important class of  non-Gaussian noise), by examining the qualitative changes of equilibrium states in  its most probable phase portraits.  We have found some peculiar bifurcation phenomena in contrast to the deterministic counterpart: (i) When  the non-Gaussianity parameter in L\'evy noise varies, there is either one, two or none backward pitchfork type bifurcations;  (ii) When a parameter in the  vector field varies, there are two or three forward pitchfork   bifurcations; (iii) The non-Gaussian L\'evy noise  clearly  leads to fundamentally more complex bifurcation scenarios, since in the special case of   Gaussian noise,  there is only one pitchfork bifurcation which is reminiscent of the deterministic situation.
\end{abstract}


\section{Introduction}
\noindent Despite the rapid development in many aspects of stochastic dynamical systems, the investigation of stochastic bifurcation is still in its infancy. A stochastic bifurcation may be defined as a qualitative change in the evolution of a stochastic dynamical system, as  a parameter varies. Stochastic bifurcations have been observed in a wide range of nonlinear systems in physical science and engineering.  The existing works on stochastic bifurcation mostly are for stochastic dynamical systems with Gaussian noise and  focus on the qualitative changes in stationary probability densities \cite{Sri1990} as solutions of steady Fokker-Planck equations,   invariant measures (together with their supports and Lyapunov  spectra) and     random point attractors \cite{Arnold2003}, or Conley index \cite{Chen2009}.



 Random  fluctuations are often assumed to have Gaussian distributions  \cite{Gui2016, Suel2006, Hasty2000b, Liu2004, Li2014} and are represented by Brownian motion.
But   the fluctuations   in some complex systems, such as temperature    evolution in   paleoclimate ice-core  records \cite{Ditlevsen1999} and  bursty transition   in gene expression \cite{Kumar2015, Dar2012},  are not Gaussian.  Then it is more appropriate to model these  random fluctuations by a non-Gaussian L\'evy motion  (i.e., $\alpha-$stable L\'evy motion)  with heavy tails and bursting sample paths \cite{Zheng2016, Klafter2011, Woyczynski2001, Chechkin2007}.

A bifurcation in deterministic low dimensional dynamical systems   often appears as a qualitative change in phase portraits in state space, and   is usually illustrated  via a bifurcation  diagram in a `parameter-steady state plane' \cite{Guckenheimer1983,Wiggins2003,Strogatz1994}.

In this present work, we   study stochastic bifurcation in a kind of stochastic phase portraits.
However, phase portraits for stochastic differential equations are delicate objects. It turns out that the phase portraits in terms of most probable orbits   \cite{Duan2015,Cheng2016} offer a promising option.  Thus  we propose here  to study   stochastic bifurcation by examining the qualitative changes (especially the changes in the number and stability type for equilibrium states)  in most probable phase portraits. To this end, we consider    bifurcation for the prototypical   scalar stochastic  differential equation with multiplicative  $\alpha-$stable L\'evy motion
$$
dX_t = (r X_t -X_t^3) dt + X_t dL_t^\alpha,
$$
where $r$ is a real parameter and the  parameter $\alpha$ is in the interval $(0, 2)$. The $\alpha-$stable L\'evy motion  $L_t^\alpha$ will be reviewed in the next section.

The deterministic counterpart $\dot x= r x-x^3$  has the well-known (forward)  `pitchfork' bifurcation \cite{Guckenheimer1983}, as the parameter $r$ increases.

\begin{figure}[!h]
\centering\includegraphics[width=0.45\textwidth]{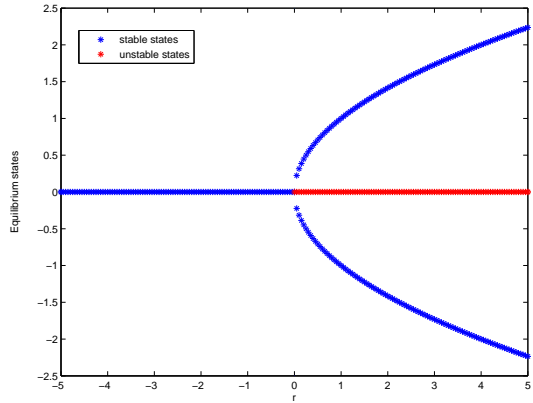}
\caption{ (Color online)  Bifurcation diagram for deterministic dynamical   system  $\dot x= r x-x^3$:  Equilibrium states vs. parameter $r$.  This is a pitchfork bifurcation at $r=0$.}
\label{Fig.0}
\end{figure}
Figure  \ref{Fig.0} is the bifurcation diagram for this deterministic  pitchfork system. For $r \leq 0$, $x=0$ is the only  equilibrium state which is stable. While for $ r>0 $, there exist two stable equilibrium states $ \sqrt r$ and $ -\sqrt r$ and one  unstable equilibrium state $x=0$.  The bifurcation parameter value is at $r=0$.

This paper is organized as follows. In Section 2, we review the definition of a scalar stable L\'evy motion $L_t^{\alpha}$ , the most probable phase portraits, and the numerical methods for   bifurcation diagrams. In Section 3, we show bifurcation diagrams for a stochastic pitchfork bifurcation under multiplicative stable L\'evy motion. Finally, we summarize our results in Section 4.

\section{Methods}
\label{methods}

\subsection{Stable L\'evy motion }
A scalar stable L\'evy motion $L_t^{\alpha}$, for $0<\alpha<2$,  is a non-Gaussian stochastic process with the following properties \cite{Duan2015,Applebaum2009,Sato1999,Samorodnitsky1994}: \\
(i)  $L_0^{\alpha} = 0$, almost surely (a.s.);\\
(ii) $L_t^{\alpha}$ has independent increments;\\
(iii)$L_t^{\alpha}$ has stationary increments: $L_t^{\alpha}-L_s^{\alpha} $ has probability distribution $S_\alpha((t-s)^\frac{1}{\alpha}, 0, 0)$ for   $  s \leq t  $; in particular,  $L_t^{\alpha}$ has distribution $S_\alpha(t^\frac{1}{\alpha}, 0, 0)$; \\
(iv) $L_t^{\alpha}$ has stochastically continuous sample paths, i.e.,   $L_t^{\alpha} \rightarrow L_s^{\alpha}$ in probability, as $t\rightarrow s$.

Here $S_{\alpha}(\sigma,\beta,\mu)$ is the so-called stable distribution  \cite{Samorodnitsky1994,Duan2015} and is  determined by four parameters, non-Gaussianity index
 $\alpha  (0 < \alpha < 2)$, skewness index $\beta (-1\leq \beta \leq 1)$, shift index $\mu  (-\infty < \mu < +\infty)$ and scale index $\sigma  (\sigma \geq 0)$.

 The stable L\'evy motion $L_t^\alpha$ has  the jump measure
 $$
\nu_{\alpha}(dy)=C_\alpha |y|^{-(1+\alpha)}\, dy,
$$
 where the coefficient
 $$
C_{\alpha} =
\frac{\alpha}{2^{1-\alpha}\sqrt{\pi}}
\frac{\Gamma(\frac{1+\alpha}{2})}{\Gamma(1-\frac{\alpha}{2})}.
$$


Note that the well-known Brownian motion $B_t$  is a special case corresponding to $\alpha=2$. Brownian motion $B_t$   has  independent and stationary increments, and has  continuous sample paths (a.s.). Moreover, $B_t -B_s $ has normal distribution $ \mathcal{N}(t-s, 0)$ for   $  s \leq t  $. In particular, $B_t$ has normal distribution $ \mathcal{N}(t, 0)$.   That is,  Brownian motion  is a Gaussian process.

\subsection{Nonlocal Fokker-Planck equation and numerical methods}

Consider a scalar stochastic differential equation with multiplicative L\'evy noise
\begin{equation} \label{sde}
  d X_t = f(X_t) dt + \sigma (X_t) d L_t^\alpha,   \;\; X_0= x_0,
\end{equation}
where $f$ is a given   vector field (or drift) and $\sigma$ is the noise intensity.

The generator for this  stochastic differential equation is
\begin{equation}
A\varphi(x)=f(x)\varphi'(x)  +  \int_{\mathbb{R}^{1}\backslash \{0\}}[\varphi(x + y\sigma(x)) - \varphi(x)] \nu_\alpha(dy).    \label{gener1}
\end{equation}
Let $ z= y\sigma(x)$.  The generator becomes
\begin{equation*}
A\varphi(x)=f(x)\varphi'(x)  + |\sigma(x)|^\alpha  \int_{\mathbb{R}^{1}\backslash \{0\}}[\varphi(x + z) - \varphi(x)]\nu_\alpha(dz).  \label{gener2}
\end{equation*}

The Fokker-Planck  equation for this stochastic differential equation, i.e., the probability
density $p(x,t)$ for the solution process $X_t $ with initial condition $X_0=x_0$ is \cite{Duan2015}
\begin{equation} \label{fpe}
 p_t  = A^* p,  \;\;   p(x,0)=\delta(x-x_0),
\end{equation}
where  $A^*$ is the adjoint operator of the generator  $A$  in   Hilbert space $ L^2(R^1) $, as defined by
\begin{equation*}
\int_{\mathbb{R}^{1}\backslash \{0\}} A\varphi(x)u(x)dx = \int_{\mathbb{R}^{1}\backslash \{0\}}\varphi(x)A^*u(x)dx.
\end{equation*}
Then via integration by parts, we get the adjoint operator for $A$
\begin{equation}
A^*u(x)=\int_{\mathbb{R}^{1}\backslash \{0\}} [|\sigma(x-y)|^\alpha u(x-y)- |\sigma(x)|^\alpha     u(x)] \; \nu_\alpha(dy).  \label{hilbert}
\end{equation}
Thus we have the nonlocal Fokker-Planck equation
\begin{equation} \label{fpe2}
p_t = - (f(x)p(x, t))_x +\int_{\mathbb{R}^{1}\backslash \{0\}} [|\sigma(x-y)|^\alpha p(x-y, t)- |\sigma(x)|^\alpha p(x, t)] \; \nu_\alpha(dy).
\end{equation}

When the stable L\'evy motion is replaced by Brownian motion, we have the following stochastic differential equation
\begin{equation} \label{sde2}
  d X_t = f(X_t) dt + \sigma (X_t) d B_t,  \;\;  X_0= x_0.
\end{equation}
The corresponding   Fokker-Planck equation is a local partial differential equation
\begin{equation} \label{fpe3}
p_t =-(f(x)p(x, t))_x + \frac12 (\sigma^2(x) p(x, t))_{xx},  \;\;   p(x,0)=\delta(x-x_0).
\end{equation}

We use a    numerical finite difference method developed  in Gao et al. \cite{Gao2016} to simulate the nonlocal Fokker-Planck  equation  (\ref{fpe2})  and use the standard finite difference method to simulate the  local Fokker-Planck  equation  (\ref{fpe3}).

\subsection{Most probable phase portraits}

As the solution of the Fokker-Planck equation,  the probability density function $p(x,t)$ is a surface in the $(x,t,p)-$space. At a given time instant $t$, the maximizer $x_m(t)$ for $p(x,t)$ indicates the most probable (i.e., maximal likely) location of this orbit at time $t$. The orbit traced out by $x_m(t)$ is called a most probable orbit  starting at $x_0$. Thus,  the deterministic orbit $x_m(t)$ follows the top ridge   of the surface in the $(x,t,p)-$space as time goes on.
For more information, see \cite{Duan2015, Cheng2016}.



\emph{Definition}: A most probable equilibrium state is a state which either attracts or repels  all nearby orbits. When it attracts    all nearby orbits, it is called a   most probable \emph{stable} equilibrium state, while if it repels all nearby orbits, it is called a   most probable \emph{unstable} equilibrium state.

A phase portrait for a stochastic dynamical system, in the sense of  most probable orbits, consists of representative orbits (including invariant objects such as most probable equilibrium states) in the state space.   Both most probable phase portraits and most probable equilibrium states are  deterministic  geometric objects.  As in the study of bifurcation for deterministic dynamical systems \cite{Guckenheimer1983, Wiggins2003, Strogatz1994}, we examine the qualitative changes in the most probable phase portraits as a parameter varies.  A simple qualitative change is the change in the `number' and `stability type' of  `most probable equilibrium states'.


\section{Results}

We now investigate the bifurcation for the scalar stochastic differential equation with multiplicative L\'evy noise
 \begin{equation} \label{pitchfork}
  d X_t = f(r, X_t) dt +  X_t \; d L_t^\alpha,
\end{equation}
 where $ f(r, X_t)= r X_t-X_t^3 $,  $r$ is a real parameter, and the non-Gaussianity parameter $\alpha \in (0, 2)$. We also compare this bifurcation diagram with that of the same system under  multiplicative Brownian noise
  \begin{equation} \label{pitchforkBM}
  d X_t = f(r, X_t) dt +  X_t \; d B_t.
\end{equation}

\emph{The existing relevant works}.  The stochastic   bifurcation for  $d X_t = f(r, X_t) dt +B_t$, with  \emph{additive} Brownian noise, was studied  in \cite{Crauel1998, Callaway2017} by examining the qualitative changes in  invariant measure and their spectral stability. The stochastic   bifurcation  for  $d X_t = f(r, X_t) dt + X_t \; B_t$, with \emph{multiplicative} Brownian noise,  was considered  in  \cite{Xu1995}  by examining the qualitative changes in  invariant measures with supports, and in \cite{Wang2015} by examining the qualitative changes in  random complete quasi-solutions. Moreover, the stochastic   bifurcation for  $d X_t = f(r, X_t) dt +L_t^\alpha$, with  \emph{additive} L\'evy noise, was studied  in \cite{ChenHQ} by considering steady probability distributions for the solutions.

\subsection{Bifurcation diagram: System under  stable L\'evy motion $L_t^\alpha$}

In the present work, we consider the case for a stochastic   bifurcation in system  (\ref{pitchfork}), with \emph{multiplicative}  $\alpha-$ stable L\'evy motion,  using  most probable phase portraits (especially most probable equilibrium states) as  a parameter $r$ in vector field  or  the non-Gaussianity parameter $\alpha$ varies. As the analytical results for most probable equilibrium states  are lacking at this time \cite{Cheng2016}, we conduct numerical simulations to generate bifurcation diagrams.

\begin{figure}[!htp]
\subfigure[]{ \label{Fig.sub.11}
\includegraphics[width=0.45\textwidth]{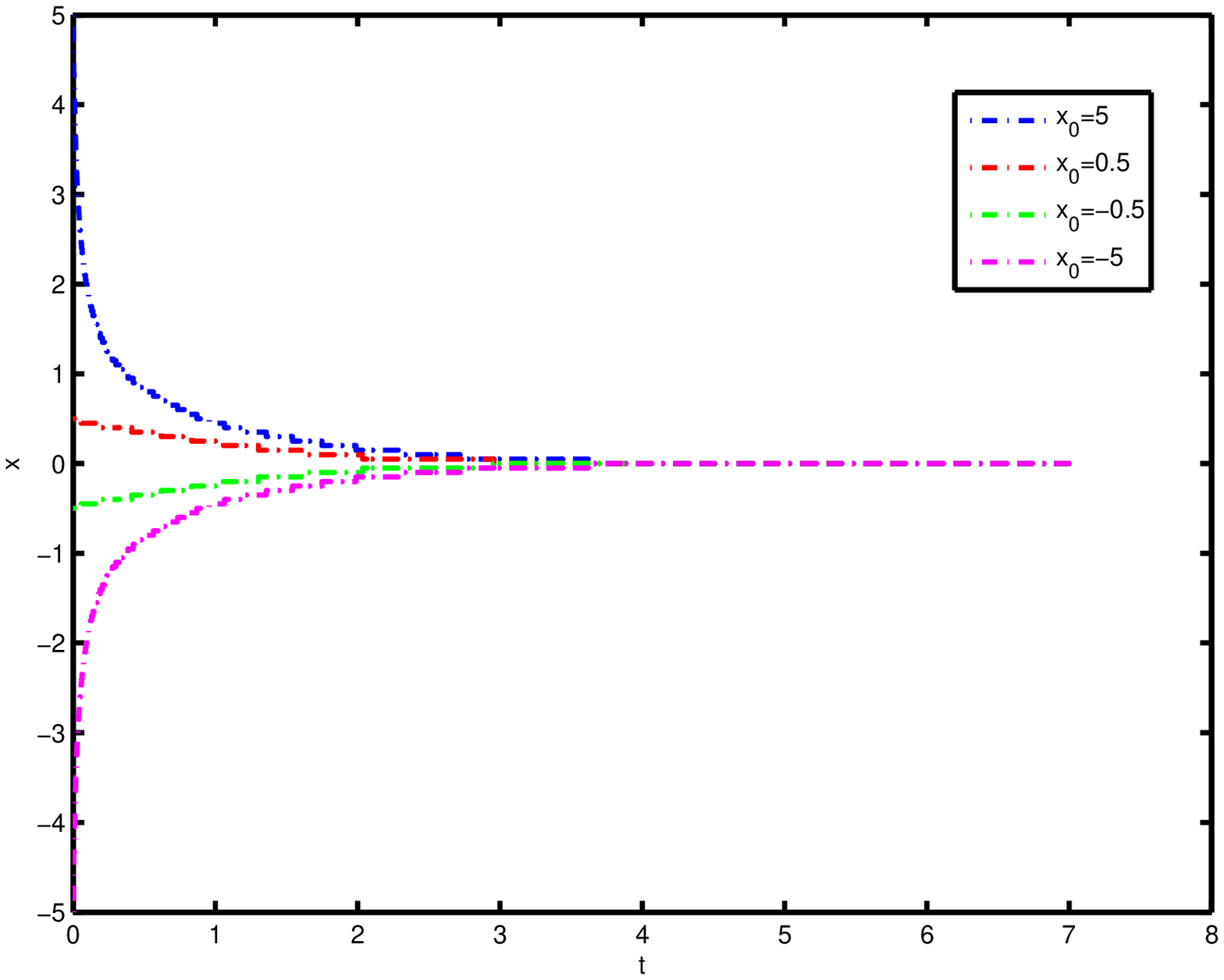}}
\subfigure[]{ \label{Fig.sub.12}
\includegraphics[width=0.45\textwidth]{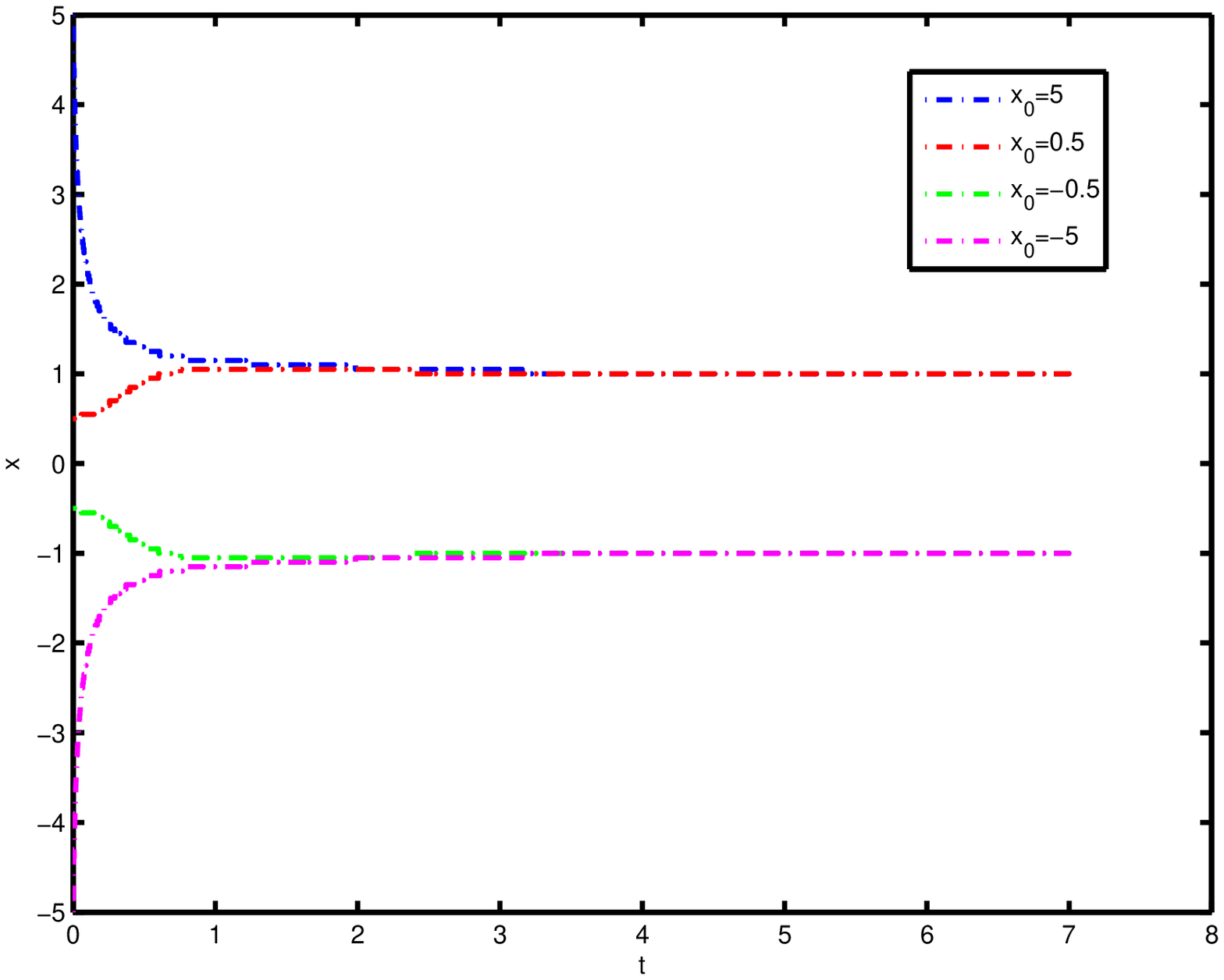}}
\caption{(Color online)  Most probable orbits and `most probable  equilibrium states' for system (\ref{pitchfork}):  (a) $ \alpha=0.3, r=-0.9 $, together with equilibrium state $x_m=0$. (b) $ \alpha=0.3, r=0.8 $, together with equilibrium states $x_m \thickapprox 1 $ and $x_m \thickapprox -1 $.}
 \label{Fig_1}
\end{figure}

For this system (\ref{pitchfork}), $0$ is always a most probable equilibrium state.
Figure \ref{Fig_1} shows the most probable orbits, starting from several initial points, with one or two most probable equilibrium states.  To generate a bifurcation diagram, we plot all possible  equilibrium states versus a parameter $r$ or $\alpha$ in the `parameter-equilibrium states plane'.

\begin{center}
\begin{figure}[!ht]
\subfigure[]{ \label{Fig.sub.21}
\includegraphics[width=0.45\textwidth]{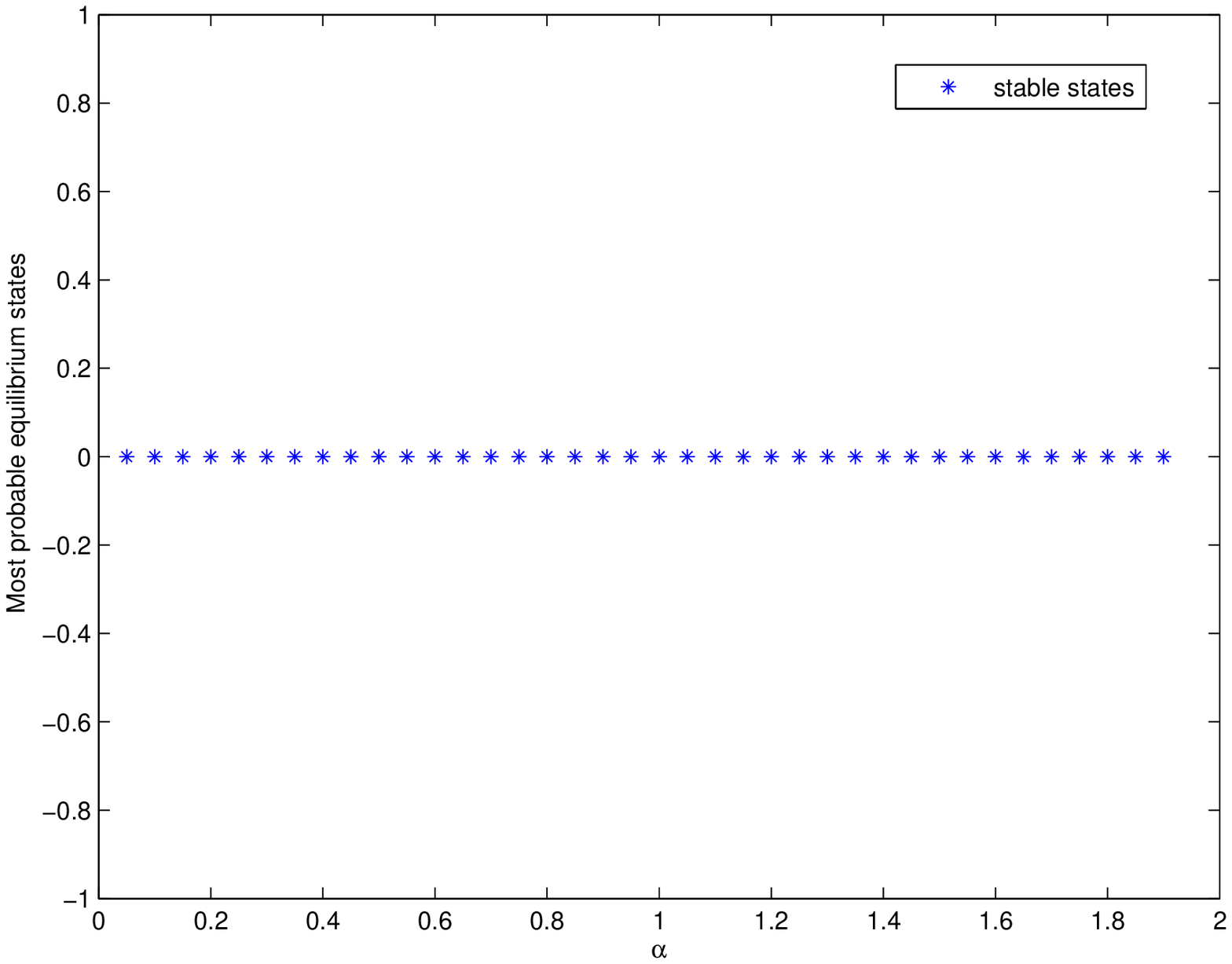}}
\subfigure[]{ \label{Fig.sub.22}
\includegraphics[width=0.45\textwidth]{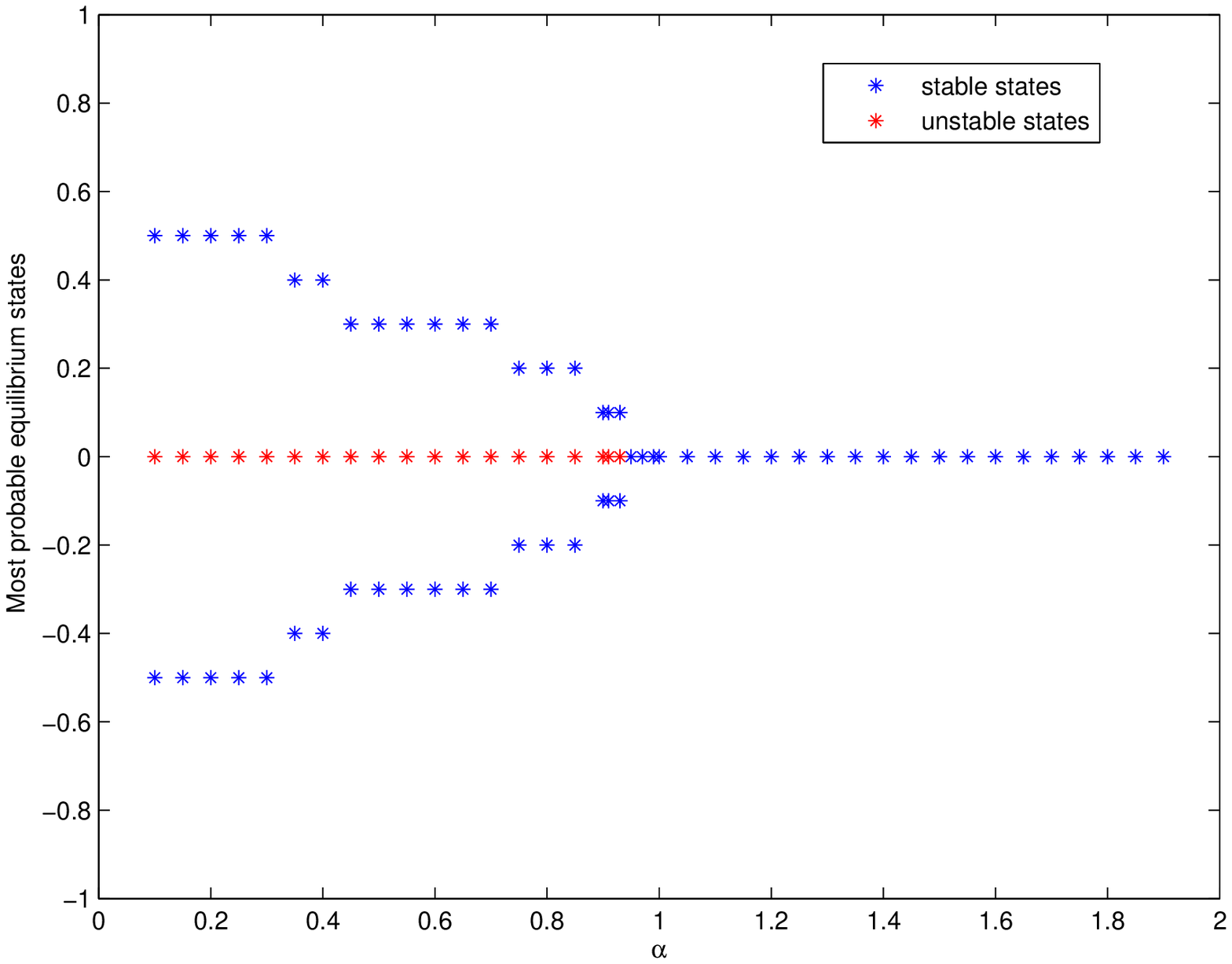}}
\subfigure[]{ \label{Fig.sub.23}
\includegraphics[width=0.45\textwidth]{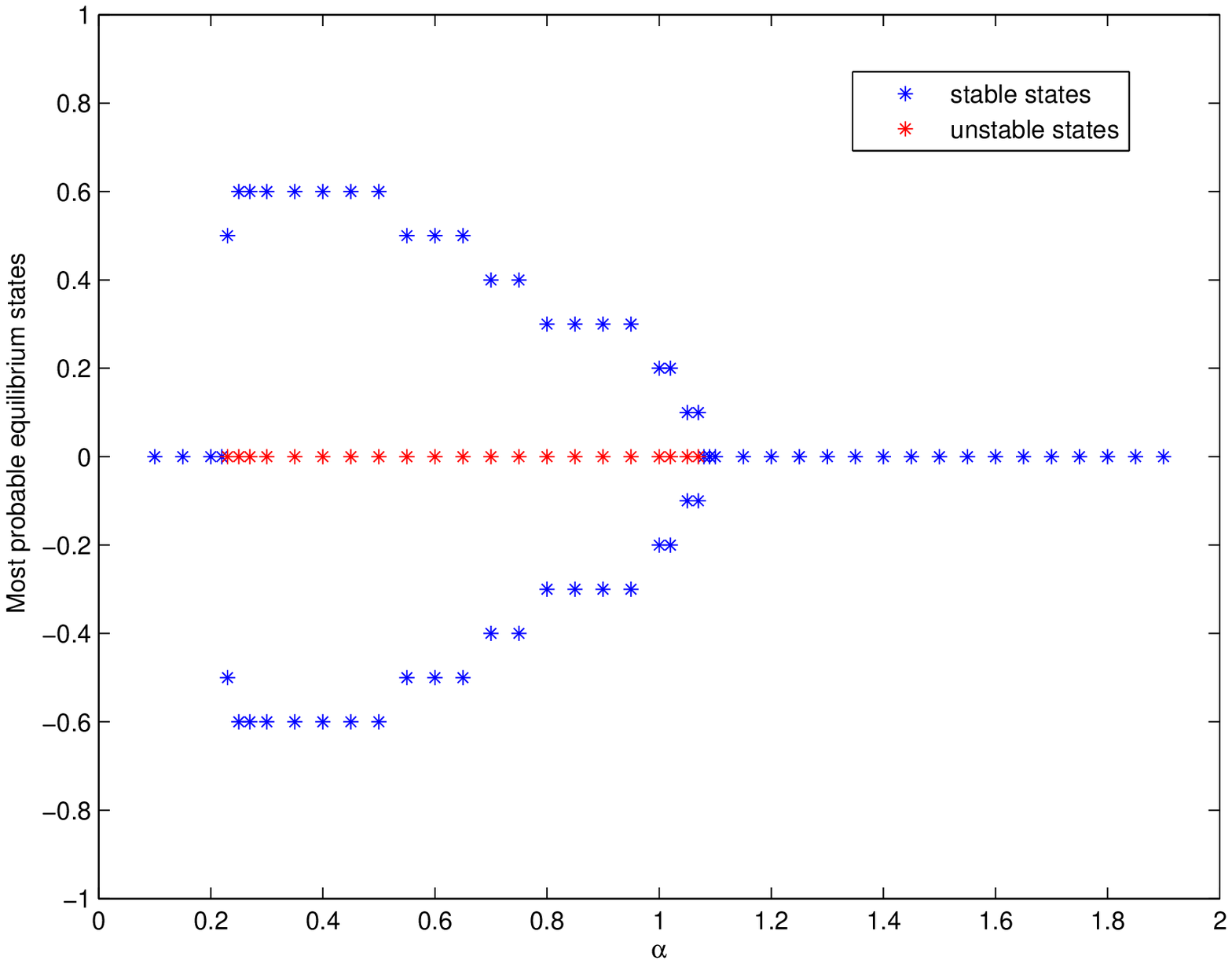}}
\subfigure[]{ \label{Fig.sub.24}
\includegraphics[width=0.45\textwidth]{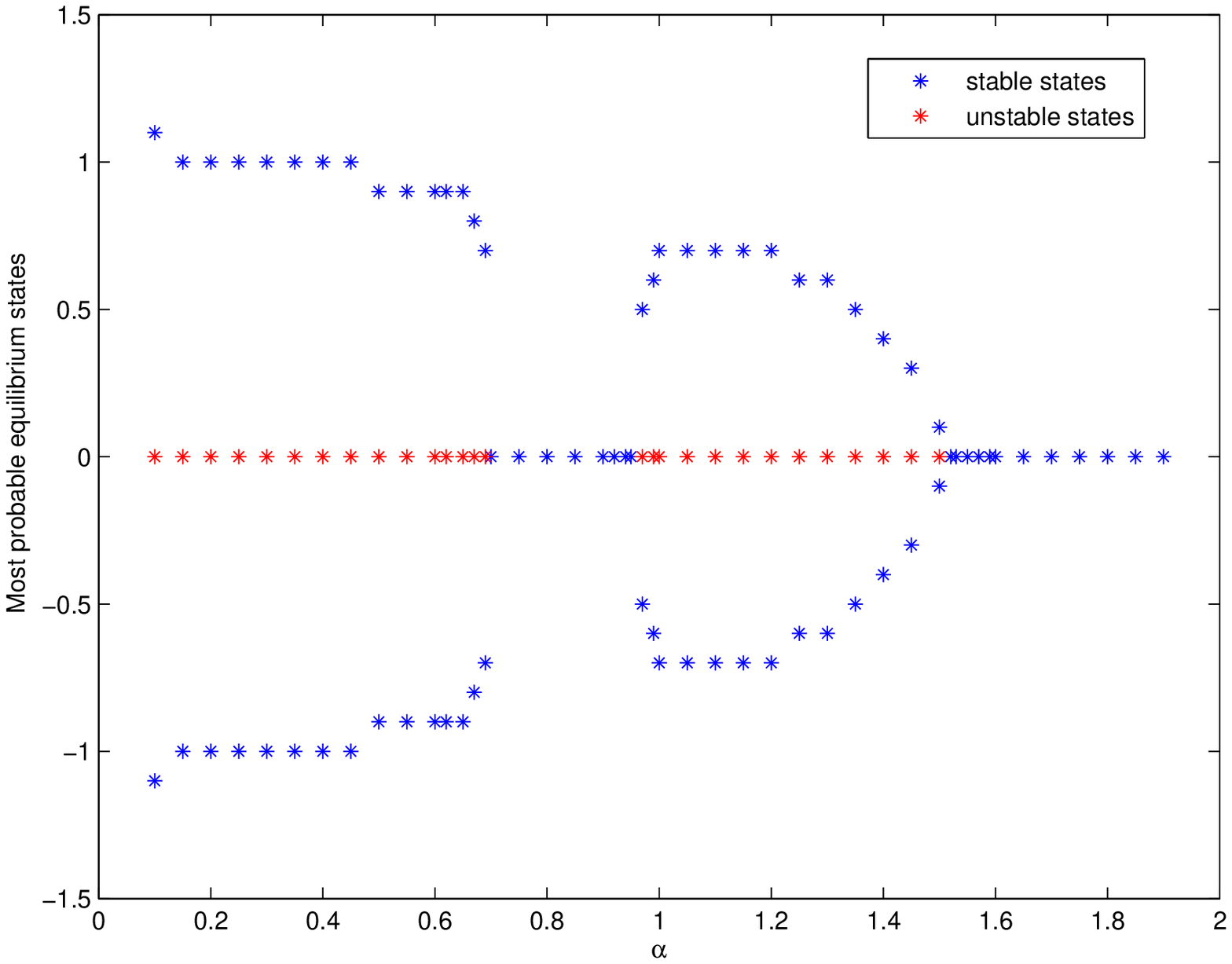}}
\subfigure[]{ \label{Fig.sub.25}
\includegraphics[width=0.45\textwidth]{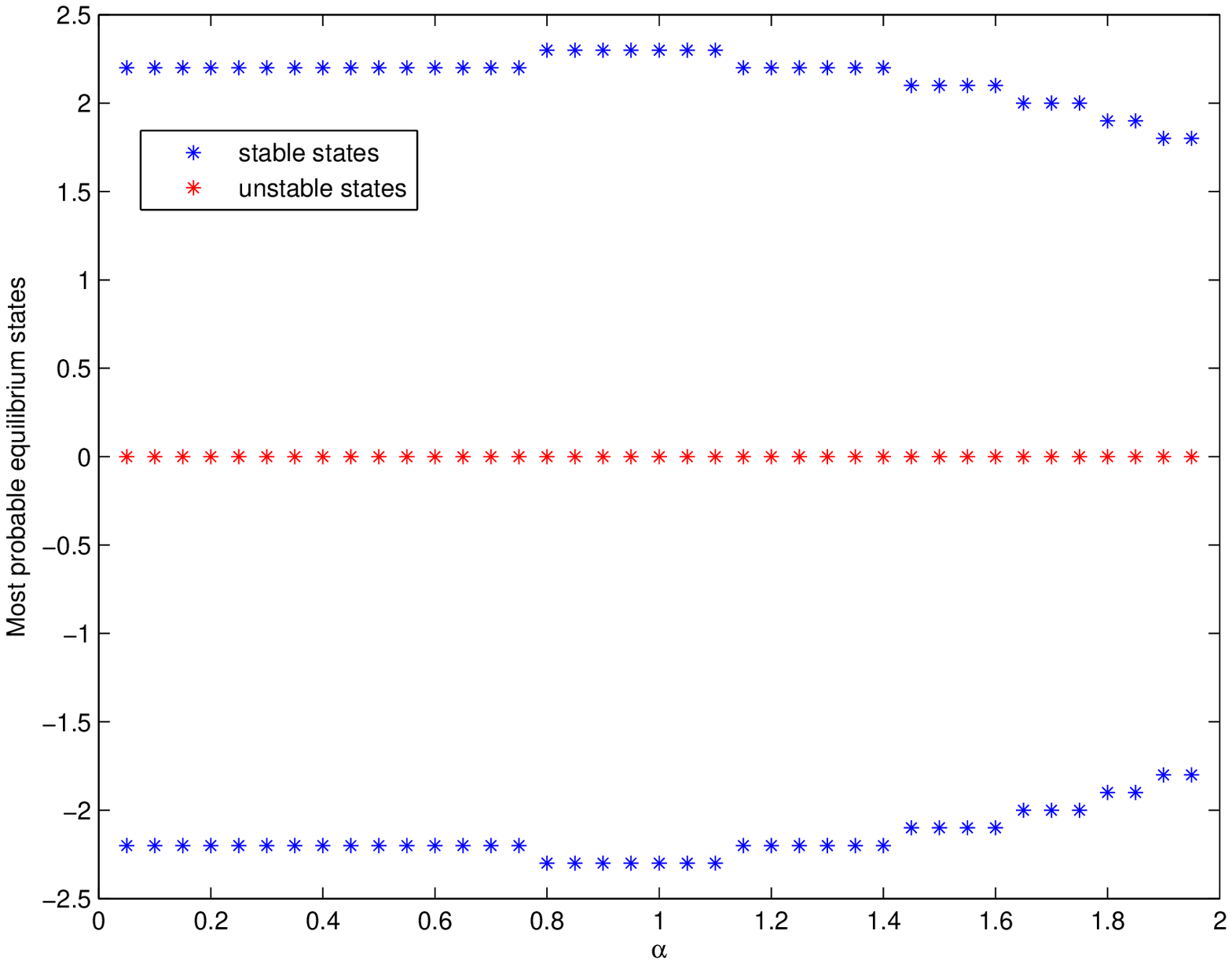}}
\caption{(Color online)  Bifurcation diagram  for system  (\ref{pitchfork}) with respect to non-Gaussianity parameter $\alpha$ :  (a) $   r  \lesssim -0.5$ (showing here $ r=-0.8 $). (b) $-0.5 \lesssim r \lesssim -0.2$ (showing here $ r=-0.2 $).  (c) $ -0.2 <  r  \lesssim 0.2$ (showing here $ r=0 $). (d) $ 0.2 \lesssim r < 2.5 $ (showing here $ r=0.8 $). (e) $ r \gtrsim 2.5  $ (showing here $ r=5 $). }
 \label{Fig_2}
\end{figure}
\end{center}

Figure \ref{Fig_2} shows the most probable equilibrium  states with respect to $\alpha$. We divide the real line $r$ into five  intervals, in each  interval the system  (\ref{pitchfork}) has the same   bifurcation phenomenon.\\
(a) For $  r  \lesssim -0.5$, the system  (\ref{pitchfork})  has only the stable equilibrium state $0$ with all $\alpha$ and there is no bifurcation.   \\
(b) For  $-0.5 \lesssim r \lesssim -0.2$,  there is a backward pitchfork bifurcation  at  $\alpha_1 \thickapprox 0.93$:  with two stable equilibrium states and one unstable equilibrium state $0$ when $\alpha < \alpha_1$ but only one unstable  equilibrium state $0$.\\
 (c) For  $ -0.2 < r  \lesssim 0.2$,  there is a   backward pitchfork bifurcation  at  $\alpha_{21} \thickapprox 1.07$  (two stable equilibrium states and one unstable equilibrium state $0$). Then there is  a `collapsing' bifurcation (as if three equilibrium states collapse into one)  at  $\alpha_{22} \thickapprox 0.23$ when two stable equilibrium states disappear but the  equilibrium state $0$ remains and becomes stable. \\
 (d)  For  $ 0.2 \lesssim r < 2.5 $,  there is a backward pitchfork bifurcation at $\alpha_{31} \thickapprox 0.69$,  a forward pitchfork bifurcation at $\alpha_{32} \thickapprox 0.95$,  and finally a `collapsing'  bifurcation at  $\alpha_{33} \thickapprox 1.5$ when two stable equilibrium states disappear but the equilibrium state $0$ remains and becomes stable. \\
(e) For $ r  \gtrsim 2.5$, the system  (\ref{pitchfork})  has two stable equilibrium states and one unstable equilibrium state $0$ and there is no bifurcation.   \\

\begin{figure}[!ht]
\subfigure[]{ \label{Fig.sub.32}
\includegraphics[width=0.45\textwidth]{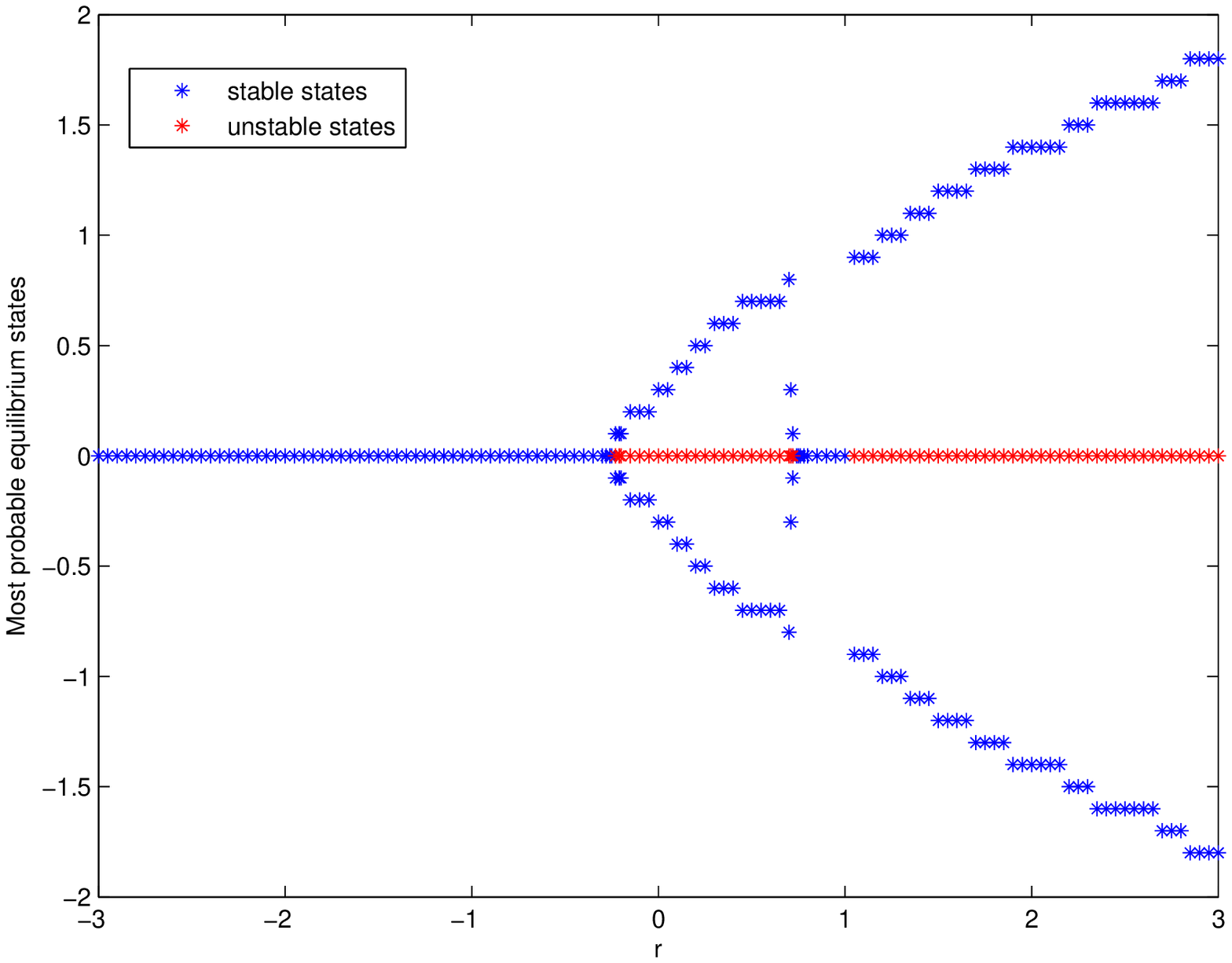}}
\subfigure[]{ \label{Fig.sub.33}
\includegraphics[width=0.45\textwidth]{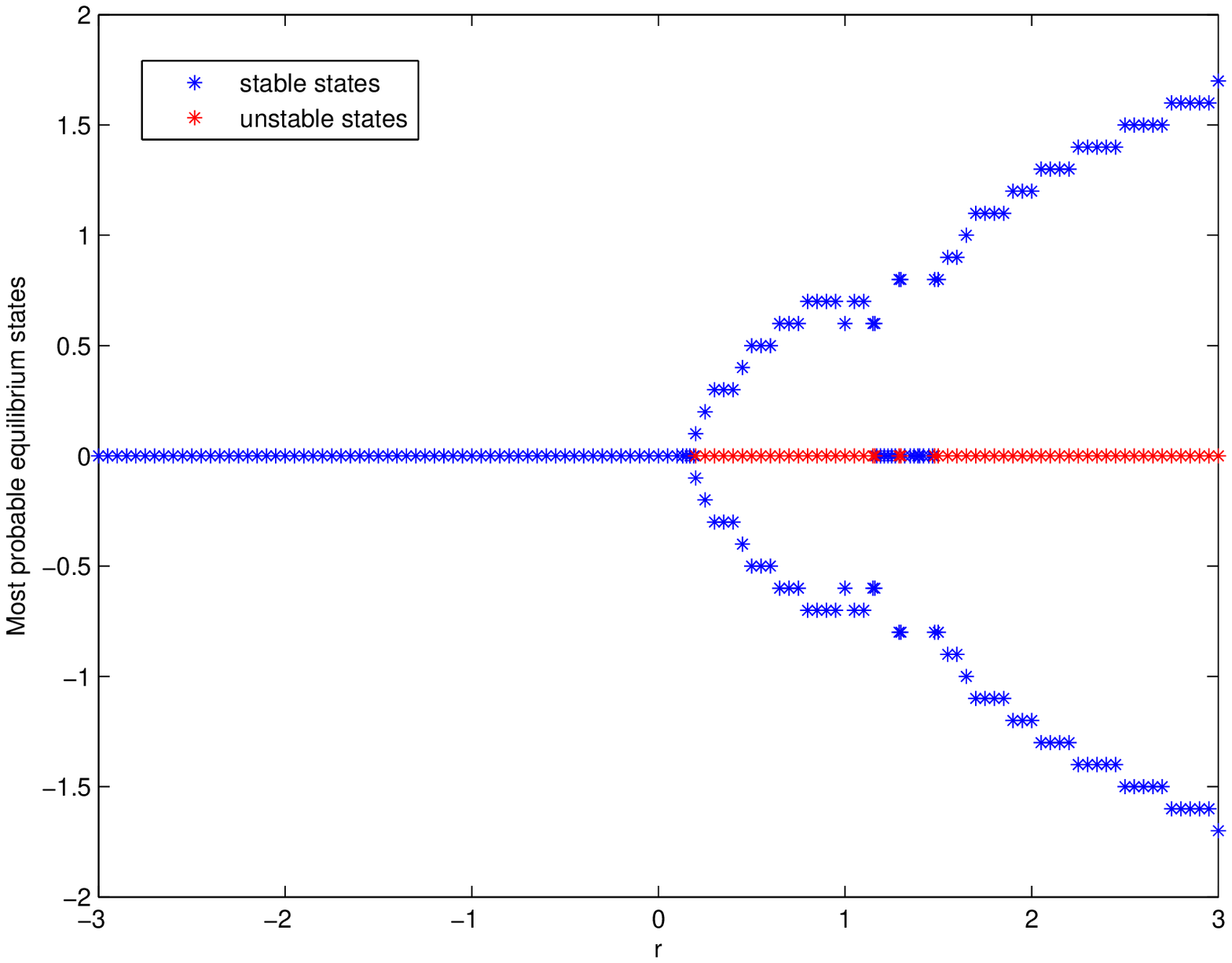}}
\caption{(Color online)  Bifurcation diagram  for system  (\ref{pitchfork}) with respect to parameter $r$ in vector field:  (a) $  0 < \alpha  \lesssim 1 $ ( showing here $ \alpha=0.9 $). (b) $1 < \alpha < 2 $  (showing here $ \alpha=1.2 $). }
 \label{Fig_3}
\end{figure}

Figure \ref{Fig_3} shows the most probable equilibrium  states with respect to $r$. The  parameter $\alpha$ can be divided into two parts, with  $\alpha=1$ as the critical or borderline value. \\
(a) For $ 0 < \alpha \lesssim 1 $, the stochastic dynamical  system  (\ref{pitchfork})   has  a forward pitchfork bifurcation at  $r_{11} \thickapprox -0.25$,  a `collapsing' bifurcation at $r_{12} \thickapprox 0.73$ when two stable equilibrium states disappear and the equilibrium state $0$ becomes stable, and finally a forward pitchfork bifurcation at $r_{13} \thickapprox 1.05$. \\
(b) For $1 < \alpha < 2 $ ,  there is a forward pitchfork bifurcation at  $r_{21} \thickapprox 0.20$ and then a `collapsing' bifurcation at  $r_{22} \thickapprox 1.17 $, suddenly a small forward pitchfork bifurcation at $r_{23} \thickapprox 1.29 $, again a `collapsing' bifurcation at $r_{24} \thickapprox 1.31 $, and finally a forward pitchfork bifurcation at $r_{25} \thickapprox 1.48 $. \\

\subsection{Bifurcation diagram: System under Brownian motion $B_t$}

\begin{figure}[!ht]
\centering\includegraphics[width=0.45\textwidth]{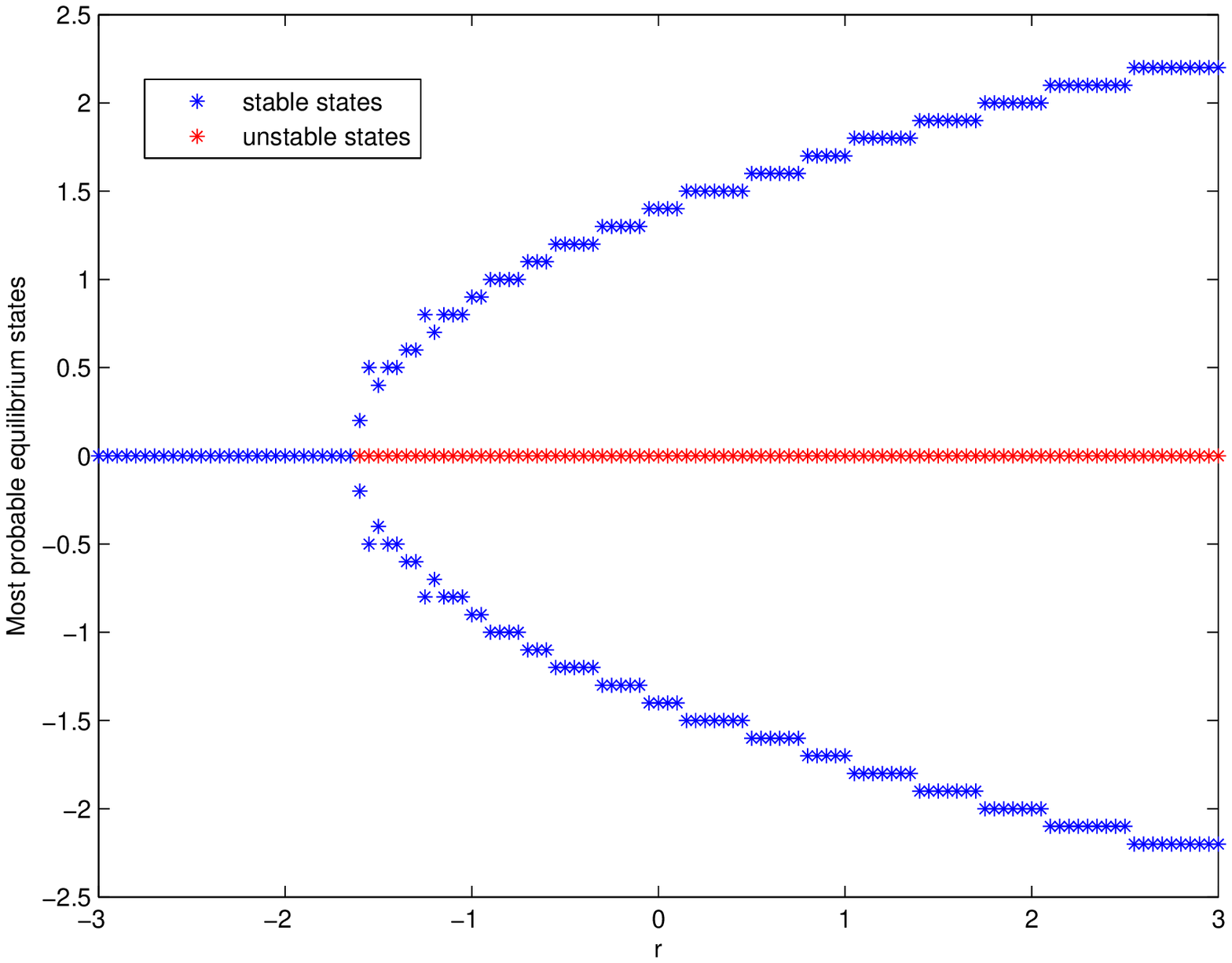}
\caption{(Color online)  Bifurcation diagram for system  (\ref{pitchforkBM})  with multiplicative Brown noise: Stochastic pitchfork bifurcation  at $r   \thickapprox -1.6$ .}
\label{Fig.4}
\end{figure}

Figure \ref{Fig.4} shows  the bifurcation diagram, i.e., the most probable equilibrium  states versus parameter $r$,   for the system  (\ref{pitchforkBM}) with multiplicative Brownian motion. There is a pitchfork bifurcation at $r \thickapprox -1.6$,  and this bifurcation was also detected in \cite[Fig. 2(b)]{Xu1995} by examining the support  of the invariant measures. This bifurcation diagram is qualitatively the same as the bifurcation diagram in Figure  \ref{Fig.0}  for the corresponding deterministic system $\dot x = rx -x^3$, although the bifurcation value is different due to the effect of noise. More significantly, this bifurcation is fundamentally different from the bifurcations under $\alpha$-stable L\'evy noise, as shown in Figure  \ref{Fig_3}.

\section{Conclusion}

 Although bifurcation studies for deterministic dynamical systems have a long history, the stochastic bifurcation investigation is still in its early stage. One reason for this slow development   in stochastic bifurcation is due to the lack of  appropriate phase portraits, in contrast to deterministic dynamical systems.

 One promising option for phase portraits of stochastic dynamical systems is the so-called most probable phase portraits \cite{Duan2015, Cheng2016}. We thus  conduct stochastic bifurcation study with the help of these phase portraits.

To demonstrate this stochastic bifurcation approach, we   study  the   bifurcation  for a system under multiplicative stable L\'evy noise (non-Gaussian). The deterministic counterpart of this system has the well-known pitchfork bifurcation.  The existing works in this topic is for the case of Brownian noise (Gaussian) and  in terms of the qualitative changes of invariant measures or point attractors.  But analytical studies of invariant measures, together with their spectra and supports, are not easily available for stochastic dynamical systems with L\'evy noise.  This also motivates us to investigate stochastic bifurcation by most probable phase portraits, especially their invariant structures such as most probable equilibrium states. By numerically examining the qualitative changes of equilibrium states in  its most probable phase portraits,  we have detected  some bifurcation phenomena such as   the double or triple  pitchfork bifurcation and a collapsing bifurcation,  when a parameter in the    vector field or in  L\'evy noise varies.

\nonumsection{Acknowledgments} \noindent
We would like to thank Xiujun Cheng and Jian Ren for helpful discussions. 
 This work was partly supported by the National Science Foundation grant 1620449, and the National Natural Science Foundation of China grants 11531006 and 11771449.


\bibliographystyle{plain}

\end{document}